\documentclass[12pt]{article}
\usepackage{amsmath,amssymb,amsbsy,amsfonts,amsthm,latexsym,
                        amsopn,amstext,amsxtra,euscript,amscd}

\begin{document}

\centerline{\sc On the logarithmic factor in error term estimates}
\centerline{\sc  in certain additive congruence problems}

\vspace{.3cm}

\centerline{{\sc M. Z. Garaev} }

\vspace{.3cm}

\begin{abstract}
We remove logarithmic factors in error term estimates in asymptotic
formulas for the number of solutions of a class of additive
congruences modulo a prime number.
\end{abstract}


\vspace{.5cm}

\centerline{\bf 1. Introduction}

\vspace{.5cm}

In additive number theory an important topic is the problem of
finding an asymptotic formula for the number of solutions of a given
congruence. In many additive congruences the error term estimates of
asymptotic formulas contain logarithmic factors. The aim of the
present paper is to illustrate application of double exponential
sums and a multidimensional smoothing argument in removing these
factors for a class of additive problems.

Let $g$ be a primitive root modulo $p$ and let $N, K$ and $M$ be any
integers with $1\le N, K <p.$  We start with recalling the well
known formula of Montgomery~\cite{Mo}
\begin{equation}
J= \frac{KN}{p} +O(p^{1/2}\log^2 p),
\end{equation}
where $J$ denotes the number of integers $x\in [H+1, H+K]$ such that
$g^x\in [M+1, M+N].$ In this paper we establish the following
statement.

\vspace{.3cm}

{\bf Theorem 1.} {\it The following formula holds:
\begin{equation}
J= \frac{KN}{p} +O\left(\frac{K^{1/2}N^{1/2}}{p^{1/4}} +
p^{1/2}\right).
\end{equation}}

\vspace{.3cm}

Formula (2) gives the asymptotic formula $J\sim KN/p$ in the range
$$
\frac{KN}{p^{3/2}}\rightarrow\infty \quad {\rm as} \quad
p\rightarrow\infty,
$$
while formula (1) gives the same asymptotic formula only when
$$
\frac{KN}{p^{3/2}\log^2 p}\rightarrow\infty \quad {\rm as} \quad
p\rightarrow\infty.
$$
Moreover, if $KN/p^{3/2}=O(1),$ then our formula guarantees the
upper bound $J\ll p^{1/2},$ while estimate (1) provides only the
bound $J\ll p^{1/2}\log^2 p.$ Thus, formula (2) includes better and
cleaner admissible range for the parameters than formula (1).

The method that we use to prove Theorem 1 is applicable to a class
of other well known additive problems. For a given integer $h,
h\not\equiv 0\pmod{p},$ denote by $J_1$ the number of solutions of
the congruence
$$
g^x-g^y\equiv h\pmod{p}, \quad 1\le x, y\le N.
$$
In \cite{Ru2} (see also \cite{Va}) the asymptotic formula
\begin{equation}
J_1 = \frac{N^2}{p} + O(p^{1/2}\log^2 p).
\end{equation}
has been established. In the present paper instead of (3) we prove

\vspace{.3cm}

{\bf Theorem 2.} {\it The following formula holds:
$$
J_1 = \frac{N^2}{p} + O\left(\frac{N}{p^{1/4}} +p^{1/2}\right).
$$}

\vspace{.3cm}

We see that in the asymptotic formula $J_1\sim N^2/p$ Theorem 2
suggests a better admissible range for $N$ than formula (3).

The following result has been obtained in \cite{Sa}:

\vspace{.3cm}

{\it Let $\mathcal{U}, \mathcal{V}\subset \{0, 1, \ldots, p-1\}$
with $u$ and $v$ elements correspondingly, and let $S$ and $T$ be
any integers with $1\le T\le p.$ If $J_2$ denotes the number of
solutions of the congruence
$$
xy\equiv z \pmod{p}, \quad x\in \mathcal{U}, \quad y\in \mathcal{V},
\quad S+1\le z\le S+T,
$$
then
\begin{equation}
\left|J_2-\frac{uvT}{p}\right|<2(puv)^{1/2}\log p.
\end{equation}}

\vspace{.3cm}

Our approach leads to

\vspace{.3cm}

{\bf Theorem 3.} {\it The following estimate holds:
$$
\left|J_2-\frac{uvT}{p}\right|\ll (u^2v^2T)^{1/3}+(puv)^{1/2}.
$$}

\vspace{.3cm}

From Theorem 3 we derive the asymptotic formula
$$
J_2\sim \frac{uvT}{p}
$$
under the condition $$\frac{uvT^2}{p^3}\rightarrow\infty \quad {\rm
as} \quad p\rightarrow\infty,$$ while estimate (4) gives the same
formula only when
$$
\frac{uvT^2}{p^3\log^2p}\rightarrow\infty \quad {\rm as} \quad
p\rightarrow\infty.
$$

We remark that estimate (4) (even with constant $2$ in the right
hand side replaced by $1$) is a consequence of the Vinogradov double
exponential sum estimate (see the Lemma below) and the inequality
$$
\sum_{a=1}^{p-1}\left|\sum_{n=S+1}^{S+T}e^{2\pi i an/p}\right|<
p^{1/2}\log p,
$$
see, for example, the proof of Lemma 5 in \cite[p. 109]{Ka2}.

\vspace{.3cm}

{\bf Theorem 4.} {\it Let $\mathcal{X}\subset \{0, 1,\ldots, p-1\}$
and let
$$
\max_{(a,p)=1}\left|\sum_{x\in \mathcal{X}}e^{\frac{2\pi
iax}{p}}\right|\le |\mathcal{X}|\Delta.
$$
If $J_3$ denotes the number of solutions of the congruence
$$
x\equiv y\pmod p, \quad x\in \mathcal{X}, \quad S+1\le y\le S+T,
$$
then
$$
J_3=\frac{|\mathcal{X}|T}{p}+O\left(\frac{|\mathcal{X}|T^{1/3}\Delta^{2/3}}{p^{1/3}}+
|\mathcal{X}|\Delta\right),
$$
where $|\mathcal{X}|$ denotes the number of elements of
$\mathcal{X}.$}

\vspace{.3cm}

Theorem 4 can be useful in estimating of discrepancy of arithmetical
sequences modulo a prime.

\vspace{.3cm}

{\bf Theorem 5.} {\it Let $h\not\equiv0\pmod p$ and let $J_4$ denote
the number of solutions of the congruence
$$
xy\equiv h\pmod p, \quad 1\le x, y\le N.
$$
Then
$$
J_4=\frac{N^2}{p}+O\left(\frac{N}{p^{1/4}}+p^{1/2}\right).
$$}

\vspace{.3cm}

In passing we remark that the argument of our paper  has also found
an application in the multiplication table problem in a residue
ring, see~\cite{GaraKara} for the details.

For more information on the latest results on distribution
properties of special sequences we refer the reader to \cite{Cobeli,
Ga, Ko, Mo, Ru2, Va} and references therein.

\vspace{.5cm}

\centerline{\bf 2. Vinogradov's double exponential sum estimate}

\vspace{.5cm}

{\sc Lemma.} {\it Let $m>1, (a, m)=1.$ Then
$$
\left| \sum_{x=0}^{m-1}\sum_{y=0}^{m-1}\nu(x)\varrho(y)e^{2\pi
i\frac{axy}{m}} \right|\le \sqrt{mXY}
$$
where $\nu(x), \varrho(y)$ are complex numbers and
$$
\sum_{x=0}^{m-1}|\nu(x)|^2=X, \quad
\sum_{x=0}^{m-1}|\varrho(y)|^2=Y.
$$}

The proof of this Lemma can be found in \cite[p. 142]{Vi}.

\vspace{.5cm}

\centerline {\bf 3. Proof of Theorem 1}

\vspace{.5cm}

If $N>p/2$ then $J$ is equal to $K$ minus the number of integers $x$
for which
$$H+1\le x\le H+K, \quad g^x\in [M+N+1, M+p]\pmod{p},$$
where now $p-N<p/2.$ For this reason it is suffice to consider the
case $N< p/2.$ By the same argument we may suppose that $K< p/2.$

Let $N_1, K_1$ be some positive integers, $ N_1<N, K_1< K$. Denote
by $J'$ the number of solutions of the congruence
$$
g^{x+z}\equiv y+t \pmod{p}
$$
subject to the conditions
$$
H+1\le x\le H+(K-K_1), \quad 1\le z\le K_1
$$
and
$$
M+1\le y\le M+(N-N_1), \quad 1\le t\le N_1.
$$
It is obvious that for fixed integers $z$ and $t$ the corresponding
number of solutions of the above congruence (in variables $x$ and
$y$) is not greater than~$J.$ Therefore,
\begin{equation}
J\ge \frac{J'}{K_1N_1}.
\end{equation}
Similarly, let $J''$ be the number of solutions to the congruence
$$
g^{x-z}\equiv y-t \pmod{p}
$$
subject to the conditions
$$
H+1\le x\le H+K+K_1, \quad 1\le z\le K_1
$$
and
$$
M+1\le y\le M+N+N_1, \quad 1\le t\le N_1.
$$
Then we have
\begin{equation}
J\le \frac{J''}{K_1N_1}.
\end{equation}

We claim that the formulas
$$
\frac{J'}{K_1N_1}= KN/p +O\left((KN)^{1/2}/p^{1/4} + p^{1/2}\right)
$$
and
$$
\frac{J''}{K_1N_1}= KN/p +O\left((KN)^{1/2}/p^{1/4} + p^{1/2}\right)
$$
hold true for some appropriately chosen parameters $K_1, N_1.$ To
prove it we express $J'$ by mean of trigonometric sums:
$$
J'=\frac{1}{p}\sum_{a=0}^{p-1}\,\sum_{x=H+1}^{H+K-K_1}\,\sum_{z=1}^{K_1}\,
\sum_{y=M+1}^{M+N-N_1}\,\sum_{t=1}^{N_1}e^{\frac{2\pi
ia(g^{x+z}-y-t)}{p}}.
$$
Picking up the term corresponding to $a=0$ we find
$$
J'=\frac{K_1 N_1 (N-N_1)(K-K_1)}{p}+ \frac{1}{p}\sum_{a=1}^{p-1}\,
\sum_{x=H+1}^{H+K-K_1}\,\sum_{z=1}^{K_1}\,\sum_{y=M+1}^{M+N-N_1}\,
\sum_{t=1}^{N_1}e^{\frac{2\pi i a(g^{x+z}-y-t)}{p}}.
$$
For $1\le a\le p-1$ we have, according to the Lemma,
$$
\left|\sum_{x=H+1}^{H+K-K_1}\,\sum_{z=1}^{K_1}e^{\frac{2\pi i
ag^{x+z}}{p}}\right|=
\left|\sum_{x=H+1}^{H+K-K_1}\,\sum_{z=1}^{K_1}e^{\frac{2\pi i
ag^xg^z}{p}}\right|\le \sqrt{pK_1(K-K_1)}.
$$
Therefore,
$$
J'-\frac{K_1 N_1 (N-N_1)(K-K_1)}{p} \ll
\frac{\sqrt{pK_1(K-K_1)}}{p}\sum_{a=0}^{p-1}
\left|\sum_{y=M+1}^{M+N-N_1}e^{\frac{2\pi i ay}{p}}\right|
\left|\sum_{t=1}^{N_1}e^{\frac{2\pi i at}{p}}\right|.
$$
The sum over $a$ is estimated by combining the Cauchy inequality
with the equalities
$$
\sum_{a=0}^{p-1}\left|\sum_{y=M+1}^{M+N-N_1}e^{2\pi
i\frac{ay}{p}}\right|^2=p(N-N_1), \quad
\sum_{a=0}^{p-1}\left|\sum_{t=1}^{N_1}e^{2\pi
i\frac{at}{p}}\right|^2=pN_1.
$$
This yields the estimate
$$
J'-\frac{K_1 N_1 (N-N_1)(K-K_1)}{p}\ll
\frac{\sqrt{pK_1(K-K_1)}}{p}\sqrt{p(N-N_1)pN_1},
$$
whence
\begin{equation}
\frac{J'}{K_1N_1}=
\frac{(N-N_1)(K-K_1)}{p}+O\left(\sqrt{\frac{pKN}{K_1N_1}}\right).
\end{equation}

If $NK < 9p^{3/2},$ then we choose $N_1=[N/2],\, K_1= [K/2]$ and
obtain
$$
\frac{J'}{K_1N_1}=O(p^{1/2})=
\frac{NK}{p}+O\left(\frac{(NK)^{1/2}}{p^{1/4}}+p^{1/2}\right).
$$

If $NK > 9p^{3/2},$ then we put
$$
\varepsilon=\frac{p^{3/4}}{(NK)^{1/2}},
$$
and observe that $\max \{1/N, 1/K\} \le \varepsilon \le 1/3.$
Therefore, in this case we can choose $N_1$ and $K_1$ such that
$$
\varepsilon N \le N_1 \le 2\varepsilon N, \quad \varepsilon K \le
K_1 \le 2\varepsilon K.
$$
Hence, from (7) we obtain
$$
\frac{J'}{K_1N_1}=\frac{NK}{p}+O\left(\frac{\varepsilon NK}{p} +
\frac{p^{1/2}}{\varepsilon}\right)=
\frac{NK}{p}+O\left(\frac{(NK)^{1/2}}{p^{1/4}}\right).
$$

Thus, in each case we have
$$
\frac{J'}{K_1N_1}=
\frac{NK}{p}+O\left(\frac{(NK)^{1/2}}{p^{1/4}}+p^{1/2}\right),
$$
whence, in view of (5), we deduce the bound
\begin{equation}
J\ge \frac{NK}{p}+O\left(\frac{(NK)^{1/2}}{p^{1/4}}+p^{1/2}\right).
\end{equation}

The above argument in application to $J''$ leads to
$$
\frac{J''}{K_1N_1}=
\frac{NK}{p}+O\left(\frac{(NK)^{1/2}}{p^{1/4}}+p^{1/2}\right),
$$
which, due to (6), yields
\begin{equation}
J\le \frac{NK}{p}+O\left(\frac{(NK)^{1/2}}{p^{1/4}}+p^{1/2}\right).
\end{equation}

The result now follows in view of (8) and (9). \qed

\vspace{.5cm}

\centerline{\bf 4. Proof of Theorem 2}

\vspace{.5cm}

We may suppose (due to (3) for example) that $N<p/2.$

Let $N_1$ be a positive integer, $N_1\le N/4.$ Denote by $J_1'$ the
number of solutions of the congruence
$$
g^{x+z}-g^y\equiv hg^{-t}\pmod{p}
$$
subject to the condition
$$
1\le x\le N-2N_1,\quad 1\le z\le N_1,\quad 1\le y\le N-N_1, \quad
1\le t\le N_1.
$$
By $J_1''$ we denote the number of solutions of the congruence
$$
g^{x-z}-g^y\equiv hg^{t}\pmod{p}
$$
subject to the condition
$$
1\le x\le N+2N_1,\quad 1\le z\le N_1,\quad 1\le y\le N+N_1, \quad
1\le t\le N_1.
$$
Then
\begin{equation}
\frac{J_1'}{N_1^2}\le J_1 \le \frac{J_1''}{N_1^2}.
\end{equation}

We express $J_1'$ in terms of trigonometric sums and then obtain
$$
J_1'=\frac{1}{p}\sum_{a=0}^{p-1}\,
\sum_{x=1}^{N-2N_1}\,\sum_{z=1}^{N_1}\,\sum_{y=1}^{N-N_1}\,
\sum_{t=1}^{N_1} e^{\frac{2\pi i a(g^{x+z}-g^y-hg^{-t})}{p}}.
$$
Using the Lemma and applying the Cauchy inequality in the same way
as in the proof of Theorem 1, we obtain
$$
\frac{J_1'}{N_1^2}= \frac{(N-2N_1)(N-N_1)}{p}+
O\left(\frac{Np^{1/2}}{N_1}\right).
$$

If $N< 10p^{3/4}$ then we let $N_1=[N/4]$ and obtain
$$
\frac{J_1'}{N_1^2}= O(p^{1/2})= N^2/p + O\left(\frac{N}{p^{1/4}}
+p^{1/2}\right).
$$

If $N > 10p^{3/4},$  then we put $N_1=[p^{3/4}]$ and obtain
$$
\frac{J_1'}{N_1^2}= \frac{N^2}{p} + O\left(\frac{N}{p^{1/2}}\right)=
N^2/p + O\left(\frac{N}{p^{1/4}} + p^{1/2}\right).
$$
Hence, for the chosen $N_1$ we have, in both cases, that
\begin{equation}
\frac{J_1'}{N_1^2}= N^2/p + O\left(\frac{N}{p^{1/4}}
+p^{1/2}\right).
\end{equation}
Analogously,
\begin{equation}
\frac{J_1''}{N_1^2}= N^2/p + O\left(\frac{N}{p^{1/4}}
+p^{1/2}\right).
\end{equation}
From (10), (11) and (12) we deduce
$$
J= N^2/p + O\left(\frac{N}{p^{1/4}} +p^{1/2}\right),
$$
whence Theorem 2. \qed

\vspace{.5cm}

\centerline{\bf 4. Proof of Theorem 3}

\vspace{.5cm}

Without loss of generality we may assume that $T\le p/2.$ Let
$T_1\le T/2$ be an integer to be chosen later. Denote by $J_2'$ the
number of solutions of the congruence
$$
xy\equiv z+t\pmod{p}
$$
subject to the condition
$$
x\in \mathcal{U}, \quad y\in \mathcal{V}, \quad S+1\le z\le S+T-T_1,
\quad 1\le t\le T_1.
$$
By $J_2''$ we denote the number of solutions of the congruence
$$
xy\equiv z-t\pmod{p}
$$
subject to the condition
$$
x\in \mathcal{U}, \quad y\in \mathcal{V}, \quad S+1\le z\le S+T+T_1,
\quad 1\le t\le T_1.
$$
Then
\begin{equation}
\frac{J_2'}{T_1}\le J_2 \le \frac{J_2''}{T_1}.
\end{equation}
Expressing $J_2'$ via trigonometric sums, picking up the main term,
applying the Lemma and the Cauchy inequality exactly in the same way
as  in the proofs of Theorems 1 and 2, we obtain
$$
\frac{J_2'}{T_1}= \frac{uvT}{p} +
O\left(\frac{uvT_1}{p}+\sqrt{\frac{uvpT}{T_1}}\right).
$$
If $uvT/p< 9 (puv)^{1/2}$ then we put $T_1=[T/2],$ and in this case
obtain
$$
\frac{J_2'}{T_1}= O\left((puv)^{1/2}\right)=\frac{uvT}{p} +
O\left((u^2v^2T)^{1/3}+(puv)^{1/2}\right).
$$
Otherwise we let $T_1$ to be any integer satisfying the condition
$$
p\left(\frac{T}{uv}\right)^{1/3} \le T_1\le
2p\left(\frac{T}{uv}\right)^{1/3}.
$$
Any such integer, due to the assumed inequality $uvT/p > 9
(uvp)^{1/2},$ also satisfies the above restriction $T_1<T/2.$ Thus,
\begin{equation}
\frac{J_2'}{T_1}= \frac{uvT}{p} +
O\left((u^2v^2T)^{1/3}+(puv)^{1/2}\right).
\end{equation}
Analogously,
\begin{equation}
\frac{J_2''}{T_1}= \frac{uvT}{p} +
O\left((u^2v^2T)^{1/3}+(puv)^{1/2}\right).
\end{equation}
Putting (13), (14) and (15) all together, we deduce Theorem 3. \qed

\vspace{.5cm}

\centerline{\bf 4. Proof of Theorem 4}

\vspace{.5cm}

We may assume that $T<p/2.$ Let $J_3'$ be the number of solutions of
the congruence
$$
x\equiv y+z\pmod p, \quad x\in \mathcal{X},\quad S+1\le y\le
S+T-T_1,\quad 1\le z\le T_1
$$
and let $J_3''$ be the number of solutions of the congruence
$$
x\equiv y-z\pmod p, \quad x\in \mathcal{X},\quad S+1\le y\le
S+T+T_1,\quad 1\le z\le T_1.
$$
Then
$$
\frac{J_3'}{T_1}\le J_3\le \frac{J_3''}{T_1}.
$$

The rest of the proof is similar to the proofs of previous theorems,
so we omit it here.

\vspace{.5cm}

\centerline{\bf 4. Proof of Theorem 5}

\vspace{.5cm} It is well known the asymptotic formula
$$
J_4=\frac{N^2}{p}+O(p^{1/2}\log^2p).
$$
Therefore, to prove Theorem 5 we may assume that $N<p/2.$

Let $J_4'$ be the number of solutions of the congruence
$$
(x+u)(y+v)\equiv h\pmod p,\quad 1\le x, y\le N-K, \quad 1\le u,v\le
K,
$$
where $K<N$ is a positive integer to be chosen later. By the same
argument that we have used in previous sections, we have the
inequality
$$
J_4>\frac{J_4'}{K^2}.
$$
Next, we express $J_4'$ in terms of trigonometric sums:
$$
J_4'=\frac{1}{p}\sum_{a=0}^{p-1}\sum_{\substack{1\le x\le N-K\\1\le
y\le N-K}}\sum_{\substack{1\le u\le K\\1\le v\le K}}e^{2\pi
i\frac{a(x+u-h(y+v)^{-1})}{p}}.
$$
Using the standard technique, we obtain
\begin{equation}
J_4'=\frac{1}{p^2}\sum_{a=0}^{p-1}\sum_{b=0}^{p-1}\sum_{z=1}^{p-1}\sum_{\substack{1\le
x\le N-K\\1\le y\le N-K}}\sum_{\substack{1\le u\le K\\1\le v\le
K}}e^{2\pi i\frac{a(x+u-hz^{-1})+b(z-y-v)}{p}}.
\end{equation}
From the classical Weil estimate of Kloosterman sums we have
\begin{equation}
\left|\sum_{z=1}^{p-1}e^{2\pi i\frac{bz-ahz^{-1}}{p}}\right|\le
2p^{1/2}
\end{equation}
for any $a\not\equiv 0\pmod p.$ If $a\equiv 0\pmod p$ and
$b\not\equiv 0\pmod p,$ then the right hand side of (17) is equal to
$-1.$ Therefore, (17) holds if at least one of the numbers $a$ and
$b$ is not divisible by $p.$ Hence, in (16) picking up the term
corresponding to $a=b=0$ and using (16) for other values of $a$ and
$b,$ we obtain
$$
J_4'=\frac{(N-K)^2K^2(p-1)}{p^2}+2\theta
p^{1/2}\left(\frac{1}{p}\sum_{a=0}^{p-1}\left|\sum_{x=1}^{N-K}e^{2\pi
i\frac{ax}{p}}\right|\left|\sum_{u=1}^{K}e^{2\pi
i\frac{au}{p}}\right|\right)^2,
$$
where $|\theta|\le 1.$ Applying the Cauchy inequality to the sum
over $a,$ we derive
$$
J_4'-\frac{(N-K)^2K^2}{p}\ll (N-K)Kp^{1/2}.
$$
Hence,
$$
J_4\ge \frac{(N-K)^2}{p}+O(NK^{-1}p^{1/2}+p^{1/2}).
$$
If $N>p^{3/4},$ then we choose $K=[p^{3/4}].$ If $N<p^{3/4},$ then
we define $K=N-1.$ In both cases we arrive at the inequality
\begin{equation}
J_4>\frac{N^2}{p}+O(Np^{-1/4}+p^{1/2}).
\end{equation}

To obtain a similar upper bound for $J_4,$ define $J_4''$ to be the
number of solutions of the congruence
$$
(x-u)(y-v)\equiv h\pmod p,\quad 1\le x, y\le N+K,\quad 1\le u,v\le
K.
$$
Then $ J_4<K^{-2}J_4''$ and
$$
J_4''=\frac{1}{p^2}\sum_{a=0}^{p-1}\sum_{b=0}^{p-1}\sum_{z=1}^{p-1}\sum_{\substack{1\le
x\le N-K\\1\le y\le N-K}}\sum_{\substack{1\le u\le K\\1\le v\le
K}}e^{2\pi i\frac{a(x-u-hz^{-1})+b(z-y+v)}{p}}.
$$
The same argument which we have used to obtain lower bounds for
$J_4'$ and $J_4$ leads to the upper bound
$$
J_4<\frac{N^2}{p}+O(Np^{-1/4}+p^{1/2}).
$$
Combining this with (18), we conclude that
$$
J_4=\frac{N^2}{p}+O(Np^{-1/4}+p^{1/2}).
$$

\vspace{.5cm}

{\sc Address of the author:

Instituto de Matem\'{a}ticas UNAM

Campus Morelia, Ap. Postal 61-3 (Xangari)

CP 58089, Morelia, Michoac\'{a}n

M\'{E}XICO}

\vspace{.5cm}

{\tt E-mail address:
 garaev@matmor.unam.mx}

\end{document}